\documentclass[11pt]{article}
\usepackage{times}
\usepackage{mathfont}
\usepackage{graphicx}
\usepackage{url}

\newtheorem{theorem}{Theorem}
\newtheorem{corollary}{Corollary}
\newtheorem{lemma}{Lemma}

\DeclareSymbolFont{lasy}{U}{lasy}{m}{n}
\let\Box\undefined
\DeclareMathSymbol\Box{0}{lasy}{"32}
\newcommand{\qed}{\hfill$\Box$}
\newenvironment{proof}{\noindent{\bf Proof:}}{\qed\medskip}

\DeclareSymbolFont{AMSb}{U}{msb}{m}{n}
\DeclareSymbolFontAlphabet{\Bbb}{AMSb}
\def\R{\ensuremath{\Bbb R}}

\let\vec\bar

\makeatletter
\long\def\@makecaption#1#2{
   \vskip 10pt 
   \setbox\@tempboxa\hbox{{\small #1. #2}}
   \ifdim \wd\@tempboxa >\hsize   
       {\small #1. #2}\par        
     \else                        
       \hbox to\hsize{\hfil\box\@tempboxa\hfil}  
   \fi}

\def\@begintheorem#1#2{\it\trivlist
  \item[\hskip\labelsep{\bf #1\ #2.\ }]}
\def\@opargbegintheorem#1#2#3{\it\trivlist
  \item[\hskip\labelsep{\bf #1\ #2\ {\rm(#3)}.}]}

\makeatother

\begin{document}

\title{Tangent Spheres and Triangle Centers}

\author{David Eppstein\thanks{Dept. of Information
\& Computer Science, Univ. of California, Irvine, CA, 92697-3425,
eppstein@ics.uci.edu}}

\maketitle

\begin{abstract}
Any four mutually tangent spheres in $\R^3$ determine three coincident
lines through opposite pairs of tangencies.
As a consequence, we define two new triangle centers.
\end{abstract}

\section{Tangent Spheres}

Any four mutually tangent spheres determine
six points of tangency.  We say that a pair of tangencies $\{t_i,t_j\}$
is {\em opposite} if the two spheres determining $t_i$ are distinct
from the two spheres determining $t_j$.  Thus the six tangencies are
naturally grouped into three opposite pairs, corresponding to the three
ways of partitioning the four spheres into two pairs.

\begin{lemma}[Altshiller-Court~{\cite[{\S}630, p.~231]{Cou-64}}]
\label{spheretan}
The three lines through opposite points of tangency of any four mutually
tangent spheres in $\R^3$ are coincident.
\end{lemma}

\begin{proof}
If three spheres have a common tangency, the three lines all meet at
that point; otherwise, each sphere either contains all of or none of the
other three spheres.
Let the four given spheres $S_i$ ($i\in\{1,2,3,4\}$)
have centers $\vec{x}_i$ and radii $r_i$.
If $S_i$ contains none of the other spheres, let $R_i=r_i^{-1}$,
else let $R_i=-r_i^{-1}$.
Then the point of tangency $t_{ij}$ between spheres $S_i$ and $S_j$
can be expressed in terms of these values as
$$t_{ij} = \frac{R_i}{R_i+R_j}\vec{x}_i + \frac{R_j}{R_i+R_j}\vec{x}_j.$$
In other words, it is a certain weighted average of the two sphere
centers, with weights inversely proportional to the (signed) radii.

Now consider the point
$$M = \frac{\sum_{i=1}^4 R_i \vec{x}_i}{\sum_{i=1}^4 R_i}$$
formed by taking a similar weighted average of all four sphere centers.
Then
$$M = \frac{R_1+R_2}{(R_1+R_2)+(R_3+R_4)}t_{12} +
\frac{R_3+R_4}{(R_1+R_2)+(R_3+R_4)}t_{34},$$
i.e., $M$ is a certain weighted average of the two tangencies $t_{12}$
and $t_{34}$, and therefore lies on the line $t_{12}\,t_{34}$.
By a symmetric argument, $M$ also lies on line $t_{13}\,t_{24}$
and line $t_{14}\,t_{23}$, so these three lines are coincident.
\end{proof}

Note that a similar weighted average for three mutually externally
tangent circles in the plane gives the Gergonne point of the triangle
formed by the circle centers.  Altshiller-Court's proof is based on the
fact that the lines
$\vec{x}_i\,t{ij}$ meet in triples at the Gergonne points of the faces of
the tetrahedron formed by the four sphere centers. We will use the
following special case of the lemma in which the four sphere centers are
coplanar:

\begin{corollary}\label{2d-opp}
The three lines through opposite points of tangency of any four mutually
tangent circles in $\R^2$ are coincident.
\end{corollary}

\section{New Triangle Centers}

\begin{figure}[t]
$$\includegraphics[height=2.5in]{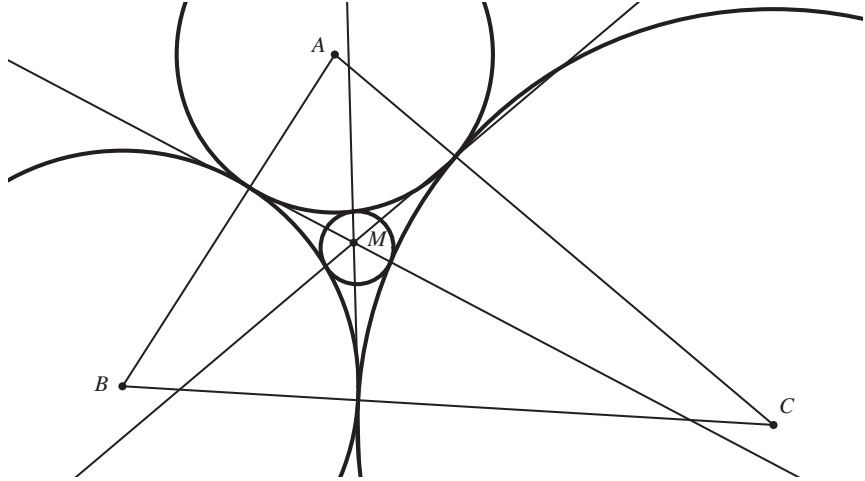}$$
\caption{A triangle $ABC$ and its new center $M$.}
\label{fig:ABC-M}
\end{figure}

Any triangle $ABC$ uniquely determines three mutually externally tangent
circles centered on the triangle vertices; if the triangle's sides have
length
$a$, $b$, $c$ then these circles have radii
$({-a+b+c})/{2}$, $({a-b+c})/{2}$, and $({a+b-c})/{2}$.
The sides of triangle $ABC$ meet its incenter at
the three points of tangency of these circles.

For any three such circles $O_A$, $O_B$, $O_C$, there exists
a unique pair of circles $O_S$ and $O_{S'}$ tangent to all three.
The quadratic relationship between the radii of the resulting two
quadruples of mutually tangent circles was famously memorialized in
Frederick Soddy's poem, ``The Kiss Precise''.

The set $\R^2\setminus(O_A\cup O_B\cup O_C)$ has five
connected components, three of which are disks and the other two of
which are three-sided regions bounded by arcs of the three circles;
we distinguish $O_S$ and $O_{S'}$ by requiring $O_S$ to lie
in the bounded three-sided region and $O_{S'}$ to lie in the unbounded
region.  Note that $O_S$ is always externally tangent to all three
circles, but $O_{S'}$ may be internally or externally tangent depending on
the positions of points $ABC$.  If $O_A$, $O_B$, and $O_C$ have a common
tangent line, then we consider $O_{S'}$ to be that line, which we think of
as an infinite-radius circle intermediate between the internally and
externally tangent cases.

We can then use Corollary~\ref{2d-opp} to define two triangle
centers: let $M$ denote the point of coincidence of the three lines
$t_{AS}\,t_{BC}$, $t_{BS}\,t_{AC}$, and $t_{CS}\,t_{AB}$
determined by the pairs of opposite tangencies of the four mutually
tangent circles
$O_A$, $O_B$, $O_C$, and $O_S$
(Figure~\ref{fig:ABC-M}), and similarly
let
$M'$ denote the point of coincidence of the three lines
$t_{AS'}\,t_{BC}$, $t_{BS'}\,t_{AC}$, and $t_{CS'}\,t_{AB}$ determined by
the pairs of opposite tangencies of the four mutually
tangent circles
$O_A$, $O_B$, $O_C$, and $O_{S'}$.
Clearly, the definitions of $M$ and $M'$ do not depend on the ordering of
the vertices nor on the scale or position of the triangle.

Despite their simplicity of definition, and despite the large amount of
work that has gone into triangle geometry~\cite{Dav-AMM-95,Kim-98},
centers $M$ and $M'$ do not appear in the lists of over 400
known triangle centers collected by Clark Kimberling and Peter Yff
(personal communications).

\section{Relations to Known Centers}

$M$ and $M'$ are not the only triangle centers defined in relation to the
``Soddy circles'' $O_S$ and $O_{S'}$.  Already known were the centers $S$
and $S'$ of these circles~\cite{Old-AMM-96,Vel-EM-66};
note that $S$ is also the point of coincidence of the three lines
$A\,t_{AS}$, $B\,t_{BS}$, $C\,t_{CS}$ and similarly for $S'$.
The Gergonne point $Ge$ can be defined in a similar way
as the point of coincidence of the three lines
$A\,t_{BC}$, $B\,t_{AC}$, and $C\,t_{AB}$.
It is known that $S$ and $S'$ are collinear with and harmonic to $Ge$ and
$I$, where $I$ denotes the incenter of triangle $ABC$~\cite{Vel-EM-66}.
Similarly $Ge$ and $I$ are collinear with and harmonic to
the isoperimetric point and the point of equal detour~\cite{Vel-AMM-85}.

\begin{theorem}
$M$ and $M'$ are collinear with and harmonic to $Ge$ and $I$.
\end{theorem}

\begin{proof}
By using ideas from our proof of Lemma~\ref{spheretan},
we can express $M$ as a weighted average of $S$ and $Ge$:
$$M = \frac{R_A A + R_B B + R_C C}{R_A+R_B+R_C+R_S} +
      \frac{R_S S}{R_A+R_B+R_C+R_S}
    = \frac{R_A+R_B+R_C}{R_A+R_B+R_C+R_S}Ge +
      \frac{R_S}{R_A+R_B+R_C+R_S}S.$$
Hence, $M$ is collinear with $S$ and $Ge$.
Collinearity with $Ge$ and $I$ follows from the known collinearity of
$S$ with $Ge$ and $I$. A symmetric argument applies to $M'$.

We omit the proof of harmonicity, which we obtained
by manipulating trilinear coordinates
of the new centers in {\em Mathematica}. See
\url{http://www.ics.uci.edu/~eppstein/junkyard/tangencies/trilinear.pdf}
for the detailed calculations.
\end{proof}

A simple compass-and-straightedge construction for the Soddy circles
and our new centers $M$ and
$M'$ can be derived from the following further relation:

\begin{theorem}
Let $\ell_A$ denote the line through point $A$, perpendicular to the
opposite side $BC$ of the triangle $ABC$.  Then the two lines
$\ell_A$ and $t_{AS}\,t_{BC}$ and the circle
$O_A$ are coincident.
\end{theorem}

\begin{proof}
Let $O_D$ be a circle
centered at $t_{BC}$, such that $O_A$ and $O_D$ cross at right angles.
Then inverting through $O_D$ produces a figure in which $O_B$ and $O_C$
have been transformed into lines parallel to $\ell_A$, while $O_A$ is
unchanged.
Since the image of $O_S$ is tangent to $O_A$ and to the two parallel
lines, it is a circle congruent to $O_A$ and centered on $\ell_A$.
Therefore, the
inverted image of $t_{AS}$ is a point $p$ where $\ell_A$ and $O_A$ cross.
Points $t_{BC}$, $t_{AS}$, and $p$ are collinear since one is
the center of an inversion swapping the other two.
\end{proof}

Since $\ell_A$, $O_A$, and $t_{BC}$ are all easy to find,
one can use this result to construct line $t_{AS}\,t_{BC}$,
and symmetrically the lines $t_{BS}\,t_{AC}$ and
$t_{CS}\,t_{AB}$, after which it is straightforward to find $O_S$, $S$,
and~$M$. A symmetric construction exists for $O_{S'}$, $S'$, and~$M'$.

\let\oldbib\thebibliography
\def\thebibliography{\addtolength{\baselineskip}{-0.5pt}\oldbib}
\bibliographystyle{abuser}
\bibliography{soddy}
\end{document}